\documentclass[12pt,leqno,A4]{article}
\usepackage{relsize}
\usepackage{paralist}
\usepackage{bibgerm}
\usepackage{makeidx,multicol}
\usepackage{makeidx}
\usepackage{epsfig}
\usepackage{graphics} 
\usepackage{textcomp} 
\usepackage{multicol}
\usepackage{amssymb} 
\usepackage{graphics} 
\usepackage{eepic}
\usepackage{empheq}
\usepackage{color} 
\usepackage{amsfonts,amsxtra}
\usepackage{amsmath,amsthm}
\usepackage{changebar}
\usepackage{enumerate}
\usepackage{nicefrac} 
\usepackage{upgreek} 
\usepackage{amscd}
\usepackage{chapterbib}
\usepackage[english]{babel}

\newcommand{\ffrac}{{\textstyle\frac{1}{2}}}
\newfont{\cfr}{eufm10}

\newcommand{\R}{{\mathbb R}}
\newcommand{\Z}{{\mathbb Z}}

\newcommand{\bc}{\mbox{\boldmath$c$}}

\newcommand{\bs}{\mbox{\boldmath$s$}}
\newcommand{\ra}{\rightarrow}

\newcommand{\oa}{\overline{a}}

\newcommand{\of}{\overline{f}}    
\newcommand{\og}{\overline{g}}

\newcommand{\cA}{{\cal A}}

\newcommand{\BBox}{\hspace{1cm}\mbox{$\Box$}}

\newcommand{\vep}{\varepsilon}

\newcommand{\ppar}{\vspace*{3mm}\noindent}
\newtheorem{prop}{Proposition}[section]
\newtheorem{theo}[prop]{Theorem}

\newtheorem{Cor}[prop]{Corollary}

\theoremstyle{definition}

\begin{document}
\title{An explicit Lyapunov function\\
for reflection symmetric\\
parabolic partial differential equations\\
on the circle\\
\vspace{2ex}
{\it \normalsize -- Dedicated to the memory of Mark Iosifovich Vishik\\
in friendship and gratitude --}
}  
\author{
\\
Bernold Fiedler* \\
Clodoaldo Grotta-Ragazzo**\\
Carlos Rocha***
\vspace*{1cm}}
\date{version of April 23, 2014}
\maketitle

\vspace*{1,5cm}

*Institut f\"ur Mathematik\\
Freie Universit\"at Berlin\\
Arnimallee 3\\ 14195 Berlin, Germany\\[4mm]
**Instituto de Matem\'atica e Estat\'\i stica\\
Universidade de S\~ao Paulo\\
05508-090 S\~ao Paulo, Brazil\\[4mm]
***Instituto Superior T\'ecnico\\
 Avenida Rovisco Pais \\1049--001 Lisboa, Portugal
\thispagestyle{empty}
\newpage
\pagestyle{plain}
\pagenumbering{arabic}
\setcounter{page}{1}

\begin{abstract}
We construct an explicit Lyapunov function for scalar parabolic reaction-advection-diffusion equations under periodic boundary conditions. 
We assume the nonlinearity is even in the advection term. We follow a method originally suggested by Matano and Zelenyak for, and limited to, 
separated boundary conditions.
\end{abstract}

\section{Introduction and main result}\label{sec1}
\numberwithin{equation}{section}
We consider real scalar semilinear parabolic partial differential equations of the form
\begin{equation}\label{eq(1.1)}
u_t=u_{xx}+f
\end{equation}
in one space dimension $0<x<1$ and with $C^1$ nonlinearities $f$. 

\ppar

Heeding Mark I. Vishik's advice "nicht zu eilen" (not to rush), we only focus on existence
versus nonexistence of Lyapunov functions in the present paper. This is but one crucial element in our ongoing
quest to clarify and classify the dynamics on the global attractors of these parabolic equations, in detail and in their
simplest scalar form. See for example \cite{FiedlerRochaWolfrum, FiRoWo12} and the references there. That entire project, in turn,
is just one modest attempt to explore a tiny part of the richness of global PDE attractors as they have
been studied, for example, in the groundbreaking and monumental work of Babin and Vishik \cite{BaVi92}
and their many followers worldwide.

\ppar

Under Dirichlet or Neumann \textit{separated boundary conditions}
\begin{equation}\label{eq(1.2)}
u=0 \quad\text{or}\quad u_{x}=0
\end{equation}
at $x=0,1$ and for nonlinearities
\begin{equation}\label{eq(1.3)}
f=f(x,u)
\end{equation}
it is well-known that there exists an explicit \textit{Lyapunov function}
\begin{equation}\label{eq(1.4)}
V(u):=\int^1_0 L(x,u,u_x)dx
\end{equation}
with the \textit{Lagrange function} integrand
\begin{equation}\label{eq(1.5)}
L(x,u,u_x):=\ffrac u^2_x-F(x,u)\,.
\end{equation}
Here $F$ is a primitive function of $f$ with respect to $u$. The Lyapunov function $V$ indeed satisfies
\begin{equation}\label{eq(1.6)}
\dot{V}=\frac{d}{dt}V(u(t,\cdot))=-\int^1_0(u_t)^2dx
\end{equation}
along any classical solution $u=u(t,x)$ of (1.1). By LaSalle's invariance principle this forces convergence to equilibria for bounded solutions and $t\ra+\infty$. Adding suitable boundary terms to the Lyapunov function $V$ the result extends to separated nonlinear boundary conditions
\begin{equation}\label{eq(1.7)}
u_x=\beta(x,u)\quad \text{at } x=0,1
\end{equation}
of Robin type. Passing from strong solutions to weak solutions similar statements remain valid and identify the semiflow (\ref{eq(1.1)}) as the $L^2$-gradient semiflow of the Lyapunov function $V$. See \cite{Henry}, \cite{Pazy} for a general background, and \cite{Mora} for the case of $C^1$-nonlinearities.

\ppar

It is a little less well-known how  \cite{Zelenyak}, and later  \cite{Matano}, extended this classical result to nonlinearities
\begin{equation}\label{eq(1.8)}
f=f(x,u,u_x)
\end{equation}
which also depend on the \textit{advection term} $u_x$, again under separated boundary conditions. For the convenience of the reader we recall the beautiful argument in the precise form of \cite{Matano} in section 2. For a suitable Lagrange function $L=L(x,u,p)$ replacing  (\ref{eq(1.5)}), the Lyapunov decay property  (\ref{eq(1.6)}) gets replaced by
\begin{equation}\label{eq(1.9)}
\dot{V}=\frac{d}{dt} V(u(t,\cdot))=-\int^1_0L_{pp}(x,u,u_x)(u_t)^2dx\,,
\end{equation}
with strict convexity of $p\mapsto L(x,u,p)$, i.e.\ with positive second partial derivative
\begin{equation}\label{eq(1.10)}
L_{pp}>0.
\end{equation}
Therefore $L_{pp}$ provides the appropriate inhomogeneous $L^2$-metric to view
 (\ref{eq(1.1)}),  (\ref{eq(1.2)}), (\ref{eq(1.8)})  as a gradient semiflow.

\ppar

Under \textit{periodic boundary conditions} $x\in S^1:=\R/\Z$, alias
\begin{equation}\label{eq(1.11)}
[u]_0^1=[u_x]^1_0=0
\end{equation}
the parabolic PDE (\ref{eq(1.1)}) retains its gradient character  (\ref{eq(1.1)}) --  (\ref{eq(1.6)}) for nonlinearities $f=f(x,u)$. The presence of advection terms $u_x$, however, is able to  produce non-equilibrium \textit{time periodic solutions} $u(t,x)$. For example consider the $SO(2)$-\textit{equivariant} case
\begin{equation}\label{eq(1.12)}
f=f(u,u_x)
\end{equation}
where $u(t,x)$ is a solution of PDE  (\ref{eq(1.1)}) iff $u(t,x+\vartheta)$ is, for any fixed rotation $\vartheta\in S^1=SO(2)$. Already \cite{AngenentFiedler} have observed that spatially nonhomogeneous \textit{rotating wave solutions}
\begin{equation}\label{eq(1.13)}
u=U(x-ct)
\end{equation}
with nonvanishing wave speeds $c\neq 0$ may then occur. Indeed this only requires nonstationary 1-periodic solutions $U$ of the traveling wave equation
\begin{equation}\label{eq(1.14)}
U''+cU'+f(U,U')=0
\end{equation}
to exist. In general, convergence to equilibria for $t\ra+\infty$ is then augmented by the possibility of convergence to rotating waves. For a specific example consider the nonlinearity
\begin{equation}\label{eq(1.14a)}
f(u,p):=\lambda u(1-u^2)-cp
\end{equation}
for $\lambda >\pi$. This amounts to viewing solutions of the cubic nonlinearity
\begin{equation}\label{eq(1.15)}
f_0(u):=\lambda u(1-u^2),
\end{equation}
known as the Chafee-Infante problem \cite{ChafeeInfante}, in coordinates which rotate at constant speed $c$ around $x\in S^1$. The nonhomogeneous equilibria $U(x)$ of the Chafee-Infante problem (\ref{eq(1.15)}) then provide nonequilibrium rotating wave solutions $U(x-ct)$ of (\ref{eq(1.14)}). Of course this argument extends to any nonlinearity $f(u,p)=f_0(u)-cp$. Other examples include nonlinearities $f=f(u,p)$ with traveling wave equations (\ref{eq(1.14)}) of Van der Pol type. For general not necessarily $SO(2)$-equivariant nonlinearities $f=f(x,u,p)$, time periodic solutions $u=u(t,x)$ may arise which are not rotating waves. Still, a Poincar\'{e}-Bendixson theorem holds which emphasizes the dichotomy between equilibria and periodic solutions for $t\ra+\infty$; see \cite{FiedlerMallet-Paret}.

\ppar

With this motivation we consider the $O(2)$-\textit{equivariant case} of PDE (\ref{eq(1.1)}) with periodic boundary conditions (\ref{eq(1.11)}) in the present paper. We therefore assume the nonlinearity $f$ to be even in $p=u_x$ to also accomodate reflections $x\mapsto -x\in S^1$ on the circle. Specifically we assume 
\begin{equation}\label{eq(1.16)}
f=f(u,p):=\of(u,\ffrac p^2)
\end{equation}
with $C^1$-nonlinearity
\begin{equation}\label{eq(1.17)}
\of=\of(u,q), \quad q=\ffrac p^2\,.
\end{equation}
Arguments based on Sturm nodal properties and zero numbers as in \cite{AngenentFiedler}, 
\cite{FiedlerMallet-Paret} then show that all rotating waves are \textit{frozen} to become equilibria, i.e.\ ``rotate'' at wave speed $c=0$. See also \cite{FiedlerRochaWolfrum}. Instead we construct an explicit Lyapunov function, in the $O(2)$-case, which forces convergence to equilibria directly by LaSalle's invariance principle. Convergence to  single equilibria, in that case, has been established by \cite{Matano} already. Those arguments  essentially excluded the alternative of rotating waves and were based on Sturm nodal  properties. They did not use the explicit Lyapunov function, which we now construct to explore the gradient flow variational character of PDE (1.1) on the circle.

\ppar

To formulate our main result, theorem 1.1 below, we assume that the $O(2)$-equivariant nonlinearity $f=\of(u,q)$ of  (\ref{eq(1.16)}),  (\ref{eq(1.17)}) 
is such that the nonautonomous ODE
\begin{align}\label{eq(1.18)}
\nonumber\frac{d}{du}q&=-\of(u,q),\\[-2mm]
\hspace*{0,01mm}\\[-2mm]
\nonumber q(u_0)&=q_0
\end{align}
possesses a global solution
\begin{equation}\label{eq(1.19)}
q(u_1)=\Psi^{u_1,u_0}(q_0)
\end{equation}
for all real $q_0$, $u_0$, $u_1$. This assumption is satisfied if $\of$ grows at most linearly in $q$:
a one-sided condition like $u \of(u,q) \leq c_1(u)+c_2(u)q$ in the relevant region $q \geq 0$
 with continuous functions $c_1, c_2$, for example, prevents
blow-up of solutions to equation (1.19) in finite ``time'' $u$; see also section 2 of
\cite{Ragazzo}. 
\ppar

We define the Lagrange function $L$, alias the integrand of the Lyapunov function $V$ in  (\ref{eq(1.4)}), as 
\begin{equation}\label{eq(1.20)}
L(u,p):=\int^p_0\int^{p_1}_0\exp(F_q(u,\ffrac p^2_2))dp_2dp_1-F(u)
\end{equation}
with the abbreviations
\begin{align}\label{eq(1.21)}
\nonumber F(u)&:=\int^u_0\of(u_1,0)\exp(F_q(u_1,0))du_1\\[-2mm]
\hspace*{0,1mm}\\[-2mm]
\nonumber F_q(u,q)&:=\int^u_0\of_q(u_1,\Psi^{u_1,u}(q))du_1\,.
\end{align}
Here $\of_q=\of_q(u_1,q_1)$ denotes the partial derivative with respect to the second argument $q_1=\Psi^{u_1,u}(q)$, and not the chain rule total derivative with respect to $q$ in $q\mapsto\of_q(u_1,\Psi^{u_1,u}(q))$.

\begin{theo}
Let $\of\in C^1$ be such that the solutions (\ref{eq(1.19)}) of ODE (\ref{eq(1.18)}) exist globally.

\ppar

Then the functional
\begin{equation}\label{eq(1.22)}
V(u):=\int^1_0 L(u,u_x)dx
\end{equation}
with the Lagrange function $L$ of (\ref{eq(1.20)}) is a Lyapunov function for the parabolic PDE (\ref{eq(1.1)}) with $O(2)$-equivariant nonlinearity $f=\of(u,\frac{1}{2}u^2_x)$ under periodic boundary conditions (\ref{eq(1.11)}). More precisely
\begin{equation}\label{eq(1.23)}
\dot{V}=\frac{d}{dt}V(u(t,\cdot))=-\int^1_0 L_{pp}(u,\ffrac(u_x)^2)(u_t)^2dx
\end{equation}
holds on classical solutions $u=u(t,x)$ of (\ref{eq(1.1)}), with strict convexity of $L(x,u,p)$ in $p$, i.e.\ with positive metric coefficient
\begin{equation}\label{eq(1.24)}
L_{pp}=\exp(F_q(u,\ffrac p^2)).
\end{equation}
\end{theo}

\bigskip
In case  $\of(u,\frac{1}{2}u^2_x)=f(u)$ is independent of $q=\frac{1}{2}p^2=\frac{1}{2}u^2_x$ we have $\of_q\equiv 0$, $F_q\equiv 0$, $L_{pp}\equiv 1$ and $L(u,p)=\frac{1}{2}p^2-F(u)$ with a primitive function $F'=f$. We therefore recover the classical Lyapunov function (\ref{eq(1.4)}) -- (\ref{eq(1.6)}) in theorem 1.1. 
\ppar

The formulation (\ref{eq(1.20)}), (\ref{eq(1.21)}) of the Lagrange function $L(u,p)$ still involves multiple integrals in terms of the evolution $\Psi^{u_1,u_0}(q_0)$ of the characteristic ODE (\ref{eq(1.18)}), (\ref{eq(1.19)}) and the nonlinearity $\of$. To eliminate some of these integrals and provide a more direct expression for $L$ we define an auxiliary function $\varphi=\varphi(u,p)$ such that 
\begin{align}\label{eq(1.25)} 
\nonumber \varphi_p(u,p)&=\Psi^{0,u}_q(\ffrac p^2)\\[-2mm]
\hspace*{0,01mm}\\[-2mm]
\nonumber\varphi(u,0)&=0.
\end{align}
Here $\Psi^{0,u}_q(q)$ denotes the partial derivative of the evolution $\Psi^{0,u}(q)$ with respect to $q$.
Of course (\ref{eq(1.25)}) amounts to simple integration,
\begin{equation}\label{eq(1.26)}
\varphi(u,p):=\int^p_0\Psi_q^{0,u}(\ffrac p^2_2)dp_2\,.
\end{equation}

\medskip
\begin{Cor}
The Lagrange function $L$ of theorem 1.1 defined in (\ref{eq(1.20)}), (\ref{eq(1.21)}) can be written equivalently as
\begin{equation}\label{eq(1.27)}
L(u,p)=p\varphi(u,p)-\Psi^{0,u}(\ffrac p^2)
\end{equation}
\end{Cor}
\bigskip

Again the trivial case $\of_q=F_q=0$, $\Psi^{0,u}(q)=q+F(u)$, $\Psi^{0,u}_q(q)=1$, $\varphi(u,p)=p$ implies $L(u,p)=p^2-\frac{1}{2}p^2-F(u)=\frac{1}{2}p^2-F(u)$.
 An  explicit construction of the Lyapunov function $V$
is also possible when $f(u,p)=a(u)+b(u)p^2/2$. Then equation 
(\ref{eq(1.18)}) is linear and can be integrated. After some computations we  can express  
 the Lagrangian $L$ of $V$ in terms of integrals of functions $a$ and $b$.
For linear $a(u)$ and constant $b$ an explicit Lyapunov function was also constructed 
in \cite{marchetti} (Proposition 5.8)  using the ideas in \cite{zel} (chapter 2).

\ppar

In section  \ref{sec2} we reproduce Matano's elegant construction of the Lagrange function for $f=f(x,u,u_x)$ and indicate where the argument fails at a technical level, as it must, under periodic boundary conditions. In section \ref{sec3} we prove theorem 1.1, based on Matano's construction. 
An alternative approach can be based
on the fact that under hypothesis (\ref{eq(1.16)})
the equilibrium equation $u_{xx}+f=0$ admits a first integral; see \cite{Ragazzo}. 
Section \ref{sec4a} proves corollary 1.2. In section \ref{sec4b}  we provide an example which shows how our Lyapunov function fails on $x\in S^1$, as it must, for nonlinearities $f(x,u,p)=f(-x,u,-p)$ which admit only a single reflection rather than full $O(2)$-equivariance $f(u,p)=f(u,-p)$ alias $f=\of(u,\frac{1}{2}p^2)$. Again this is due to the occurrence of nonstationary time periodic orbits. Section \ref{sec5}  collects comments on the associated PDE global attactors, on quasilinear equations, and on negative $q=\frac{1}{2}u^2_x$ alias imaginary $u_x$.

\ppar

\textbf{Acknowledgement.} The author Bernold Fiedler is much indebted to Hiroshi Matano for drawing his attention to, 
and patiently explaining, his construction of a Lagrange function for nonlinearities $f=f(x,u,u_x)$ 
under separated boundary conditions many years ago. 
Alexey Alimov and Pavel Gurevich have significantly improved the paper, by translation.
We are also grateful for patient typesetting by M.\ Barrett.

\ppar

Bernold Fiedler and Carlos Rocha were partially supported by the Deut\-sche Forschungsgemeinschaft, SFB 647 Space -- Time -- Matter and by FCT Portugal.
Clodoaldo Grotta-Ragazzo was partially supported by CNPq Grant 305089/2009-9, Brazil.

\section{Matano's construction}\label{sec2}

In this section we recall Matano's construction \cite{Matano} of a Lagrange function $L=L(x,u,p)$ such that
\begin{equation}\label{eq(2.1)}
V(u):=\int^1_0 L(x,u,u_x)dx
\end{equation}
becomes a Lyapunov function for PDE (\ref{eq(1.1)})  under separated Dirichlet or Neumann boundary conditions (\ref{eq(1.2)}) and  for general nonlinearities $f=f(x,u,u_x)$. We show 
\begin{equation}\label{eq(2.2)}
\dot{V}=\frac{d}{dt}V(u(t,\cdot))=-\int^1_0L_{pp}(x,u,u_x)\cdot(u_t)^2dx
\end{equation}
with the strict convexity condition
\begin{equation}\label{eq(2.3)}
L_{pp}>0.
\end{equation}

\ppar

See \cite{Zelenyak}, Lemma 1 for the closely related original construction due to Zelenyak. Twenty years later the original construction was retrieved from oblivion, clarified, and slightly modified under less restrictive regularity assumptions in the appendix to \cite{Matano}.  The differences become more apparent in the more general quasilinear parabolic case (\ref{eq(5.1)}) which we review briefly in the discussion of section 6. For clarity of presentation we follow Matano here, rather than the somewhat convoluted original argument by Zelenyak. Again we emphasize that both constructions
of a Lyapunov function are limited to separated boundary conditions,
albeit of slightly more general form, and must fail for periodic boundary
conditions.

\ppar

The construction proceeds as follows. For classical solutions $u=u(t,x)$ we integrate (\ref{eq(2.1)}) by parts and substitute $u_{xx}=u_t-f$ to obtain from (\ref{eq(1.1)}) 
\begin{align}\label{eq(2.4)}
\nonumber\dot{V}&=\int^1_0(L_uu_t+L_pu_{tx})dx=\\
\nonumber&=\int^1_0(L_u-\frac{d}{dx}L_p(x,u,u_x))u_tdx=\\
&=\int(L_u-L_{px}-L_{pu}u_x-L_{pp}u_{xx})u_tdx=\\
\nonumber&=\int^1_0((L_u-L_{px}-L_{pu}u_x+L_{pp}f)u_t-L_{pp}(u_t)^2)dx=\\
\nonumber&=-\int^1_0L_{pp}(u_t)^2dx,
\end{align}
as required. Here we have assumed Dirichlet boundary conditions $u=0$, and hence $u_t=0$, at $x=0,1$, for simplicity. Neumann boundary conditions or more general nonlinear boundary conditions (\ref{eq(1.7)}) of Robin type can be covered by adding suitable boundary terms to $V$. To satisfy the last equality, of course, the Lagrange function $L$ is required to satisfy the linear first order PDE
\begin{equation}\label{eq(2.5)}
L_u-L_{xp}-pL_{up}+fL_{pp}=0
\end{equation}
for all real arguments $u$, $p$, and $0\leq x\leq 1$. To reduce the order and guarantee convexity condition $L_{pp}>0$, Matano makes the Ansatz 
\begin{equation}\label{eq(2.6)}
L_{pp}=:\exp g.
\end{equation}
Differentiating (\ref{eq(2.5)}) partially with respect to $p$, the terms $L_{up}$ cancel and he obtains the first order linear PDE
\begin{equation}\label{eq(2.7)}
g_x+pg_u-fg_p=f_p
\end{equation}
for $g=g(x,u,p)$. This linear first order PDE for $g$ can be solved by the method of characteristics: along the solutions 
$(u,p)(x)$ of the ODE
\begin{align}\label{eq(2.8)}
\nonumber \frac{du}{dx}&=p\\[-2mm]
\hspace*{0,01mm}\\[-2mm]
\nonumber\frac{dp}{dx}&=-f(x,u,p)\,,
\end{align}
the function $x \mapsto g=g(x,u(x),p(x))$ must have total derivative
\begin{equation}\label{eq(2.9)}
\frac{d}{dx}g=f_p(x,u,p).
\end{equation}
For example we may assume
\begin{equation}\label{eq(2.10)}
g(0,u,p)\equiv 0
\end{equation}
and obtain $g$ globally, in this way, provided that the solutions of the characteristic ODE (\ref{eq(2.8)}) exists for all $0\leq x\leq 1$ and for all real initial conditions $u,p$ at $x=0$. In shorthand, ascending from (\ref{eq(2.7)})$\cdot L_{pp}=(\ref{eq(2.5)})_p$ to (\ref{eq(2.5)}) itself, via (\ref{eq(2.6)}), we then define
\begin{align}\label{eq(2.11)}
\nonumber L(x,u,p)&:=\int^p_0\int^{p_1}_0\exp g(x,u,p_2)dp_2dp_1-F(x,u), \\
\hspace*{0,01mm}\\[-2mm]
\nonumber F(x,u)&:=\int^u_0 f(x,u_1,0)\exp(g(x,u_1,0))du_1\,.
\end{align}
Indeed the left hand side of (\ref{eq(2.5)}) is independent of $p$, by this construction. Therefore (\ref{eq(2.5)}) holds, for all $p$, if we verify that (\ref{eq(2.5)}) holds at $p=0$. At $p=0$, definition (2.11) implies $0\equiv L_p\equiv L_{px}$ and $L_u=-F_u=-f\exp g=-fL_{pp}$. This proves (\ref{eq(2.5)}) and completes the Matano construction of the Lyapunov function $V$.

\ppar

Of course this correct construction must fail when abused to cover periodic boundary conditions. And it does. Suppose the characteristic equation (\ref{eq(2.8)}) possesses a periodic orbit $(u,p)(x)$ of period one, i.e.
\begin{equation}\label{eq(2.12)}
[(u,p)(x)]^1_0=0.
\end{equation}
Then 1-periodicity of $x\mapsto g(x,u,p)$ requires
\begin{align}\label{eq(2.13)}
\nonumber 0&=[g(x,u(0),p(0))]^1_0=[g(x,u(x),p(x))]^1_0=\\[-2mm]
\hspace*{0,01mm}\\[-2mm]
\nonumber &=\int^1_0\frac{d}{dx}g(x,u(x),p(x))dx=\int^1_0f_p(x,u(x),p(x))dx
\end{align}
 in view of (\ref{eq(2.9)}). But this integrability condition for $f_p$ may easily be violated, keeping the periodic orbit $(u,p)(x)$ unaffected. The Matano construction must therefore fail, in general, whenever time periodic orbits appear in PDE (\ref{eq(1.1)}). The nonlocality of the ill-posed compatibility condition
\begin{equation}\label{eq(2.14)}
[g(x,u(x),p(x))]^1_0=\int^1_0f_p(x,u(x),p(x))dx
\end{equation}
along the characteristics (\ref{eq(2.8)}), however, comes to the rescue of the Matano construction in the $O(2)$-equivariant case.

\section{Proof of theorem 1.1} \label{sec3}

The proof of theorem 1.1 consists of a slight adaptation of the Matano construction, from section 2, to the case of $O(2)$-equivariant nonlinearities
\begin{equation}\label{eq(3.1)}
f=f(u,p)=\of(u,q), \quad q:=\ffrac p^2.
\end{equation}
The nonlinearity $f$ is even in $p=u_x$, due to reflections,  and does not depend on $x$ explicitly, due to rotations. Consequently we consider $x$-independent Lagrange functions $L=L(u,p)$ and seek $O(2)$-invariant Lyapunov  functions of the form
\begin{equation}\label{eq(3.2)}
V(u):=\int^1_0L(u,u_x)dx,
\end{equation}
where
\begin{equation}\label{eq(3.3)}
L_{pp}=\exp(g)>0
\end{equation}
and $g$ takes the reflection symmetric form
\begin{equation}\label{eq(3.4)}
g=g(u,p)=\og(u,q),\quad q:=\ffrac p^2\,.
\end{equation}
The Matano calculation (\ref{eq(2.4)}) -- (\ref{eq(2.7)}) then leads to the first order linear PDE
\begin{equation}\label{eq(3.5)}
p(\og_u -\of\og_q)=p\of_q\,.
\end{equation}
Here we have used the chain rule and we substituted definitions (\ref{eq(3.1)}), (\ref{eq(3.4)}) of $f,g$ in (\ref{eq(2.7)}). We divide by $p$ and solve
\begin{equation}\label{eq(3.6)}
\og_u-\of\og_q=\of_q
\end{equation}
by the method of characteristics along the global solutions $q(u_1)=\Psi^{u_1,u_0}(q_0)$ of 
\begin{align}\label{eq(3.7)}
\nonumber\frac{d}{du}q&=-\of(u,q)\\[-2mm]
\hspace*{0,01mm}\\[-2mm]
\nonumber q(u_0)&=q_0
\end{align}
defined in  (\ref{eq(1.18)}), (\ref{eq(1.19)}). Then any $\og$ which satisfies
\begin{equation}\label{eq(3.8)}
\frac{d}{du}\og(u,q(u))=\of_q(u,q(u))
\end{equation}
along the characteristics, say with initial condition
\begin{equation}\label{eq(3.9)}
\og(0,q):=0,
\end{equation}
solves the first order linear PDE (\ref{eq(3.5)}). With the help of the evolution $q(u_1)=\Psi^{u_1,u_0}(q_0)$ this implies
\begin{equation}\label{eq(3.10)}
\og(u,q)=\int^u_0\of_q(u_1,\Psi^{u_1,u}(q))du_1=:F_q(u,q);
\end{equation}
see also the abbreviation (\ref{eq(1.21)}).

\ppar

Again we ascend from (\ref{eq(3.10)}) to (\ref{eq(2.5)}) which now reads
\begin{equation}\label{eq(3.11)}
L_u-pL_{up}+fL_{pp}=0.
\end{equation}
With the definition $L_{pp}:=\exp\og=\exp F_q$ we obtain
\begin{equation}\label{eq(3.12)}
L(u,p):=\int^p_0\int^{p_1}_0(\exp F_q(u,\ffrac p^2_2))dp_2dp_1-F(u).
\end{equation}
Here $F(u)$ is a suitable integration constant. To determine $F(u)$ we only have to evaluate  (\ref{eq(3.11)}) at $p=0$ to obtain
\begin{equation}\label{eq(3.13)}
\frac{d}{du}F(u)=-L_u=fL_{pp}=f(u,0)\exp F_q(u,0).
\end{equation}

\ppar

The requirements  (\ref{eq(3.12)}), (\ref{eq(3.13)}) are satisfied for the Lagrange function $L$ defined in  (\ref{eq(1.20)}).  (\ref{eq(1.21)}). This proves theorem 1.1. \BBox

\ppar

\section{Proof of Corollary 1.2}\label{sec4a}

The proof of corollary 1.2 proceeds by performing the integrals in the definition (1.20), (1.21) of the Lagrange function $L(u,p)$, explicitly. Alternatively, of course, it is possible to verify directly that $V(u)=\int_0^1L(u,u_x)dx$ is a Lyapunov function. This would not motivate the construction of $L$, however, in contrast to Matano's elegant approach.

\ppar

To evaluate the integrals  (\ref{eq(1.20)}),  (\ref{eq(1.21)}) we first observe that the derivative $\eta(u_1)$ of the evolution $\Psi^{u_1,u_0}(q_0)$ of the characteristic ODE  (\ref{eq(1.18)}) with respect to the initial condition $q_0$, 
\begin{equation}\label{eq(4a.1)}
\eta(u_1):=\Psi^{u_1,u_0}_q (q_0),
\end{equation}
satisfies the linearized characteristic equation
\begin{align}\label{eq(4a.2)}
\nonumber \frac{d}{du_1}\eta(u_1)&=-\of_q(u_1,\Psi^{u_1,u_0}(q_0))\eta(u_1)\\[-2mm]
\hspace*{0,01mm}\\[-2mm]
\nonumber \eta(u_0)&=1.
\end{align}
Explicit integration of  (\ref{eq(4a.2)}) shows 
\begin{equation}\label{eq(4a.3)} 
\eta(u_1)=\exp\left(-\int^{u_1}_{u_0}\of_q(u_2,\Psi^{u_2,u_0}(q_0))du_2\right).
\end{equation}
Inserting $u_0=u$, $u_1=0$, $q_0=q$ we obtain
\begin{equation}\label{eq(4a.4)}
\exp F_q(u,q)=\eta(0)=\Psi^{0,u}_q(q)
\end{equation}
from  (\ref{eq(1.21)}). Insertion of  (\ref{eq(4a.4)}) with $q=\frac{1}{2}p^2_2$ in the double integral $L_1$ of definition (1.21) of the Lagrange function $L(u,p)$, and subsequent integration by parts, then yields

\begin{align}\label{eq(4a.5)}
\nonumber L_1(u,p):=&\int^p_0\int^{p_1}_0\exp(F_q(u,\ffrac p^2_2))dp_2dp_1=\\
\nonumber =&\int^{p}_0 1\cdot \int^{p_1}_0\Psi^{0,u}_q(\ffrac p^2_2)dp_2\ dp_1=\\
 =&\left[p_1\int^{p_1}_0\Psi^{0,u}_q(\ffrac p_2^2)dp_2\right]^p_0
 -\int^p_0p_1\Psi^{0,u}_q(\ffrac p_1^2)dp_1=\\
\nonumber=&\ p\varphi(u,p)-\int^{\ffrac p^2}_0\Psi^{0,u}_q(q)dq=\\
\nonumber=&\ p\varphi(u,p)-\Psi^{0,u}(\ffrac p^2)+\Psi^{0,u}(0).
\end{align}
Here we have used definition (1.26) of the auxiliary function $\varphi$ and we  have substituted $q=\frac{1}{2}p^2$ in the integral.

\ppar

To evaluate the remaining term $-F(u)$ in $L=L_1-F$ of (\ref{eq(1.20)}) we first observe that the evolution property of $\Psi^{u_1,u_0}$ trivially implies that the partial derivative $\Psi^{u_1,u}_u(q)$ satisfies
\begin{align}\label{eq(4a.6)}
\nonumber\Psi^{u_2,u}_u(q)&=\lim_{\vep\ra0}\vep^{-1}(\Psi^{u_2,u+\vep}(q)-\Psi^{u_2,u+\vep}(\Psi^{u+\vep,u}(q))\\[-2mm]
\hspace*{0,01mm}\\[-2mm]
\nonumber&= \Psi^{u_2,u}_q(q)\of(u,q),
\end{align}
by the chain rule and definition (1.18), (1.19) of the evolution $\Psi$. Therefore (\ref{eq(1.21)}), (\ref{eq(4a.4)}) and (\ref{eq(4a.6)}) with $u=u_1$, $u_2=q=0$ imply
\begin{align}\label{eq(4a.7)}
\nonumber -F(u)&=-\int^u_0\of(u_1,0)\exp(F_q(u_1,0))du_1=\\
\nonumber &=-\int^u_0\of(u_1,0)\Psi^{0,u_1}_q(0)du_1=-\int^u_0\Psi^{0,u_1}_{u_1}(0)du_1=\\   
&=-[\Psi^{0,u_1}(0)]^u_0=-\Psi^{0,u}(0).
\end{align}
Addition of  (\ref{eq(4a.5)}) and (\ref{eq(4a.7)}) implies
\[L(u,p)=L_1(u,p)-F(u)=p\varphi(u,p)-\Psi^{0,u}(\ffrac p^2)\]
as claimed in (\ref{eq(1.27)}). This proves the corollary.

\section{Reflection symmetry}\label{sec4b}
In this section we study the parabolic PDE (\ref{eq(1.1)}) under periodic boundary conditions $x\in S^1=\R/2\pi\Z$, as in (\ref{eq(1.11)}). We consider nonlinearities $f=f(x,u,u_x)$ which are required to possess only a single reflection symmetry $x\mapsto -x$, i.e.\ 
\begin{equation}\label{eq(4.1)}
f(-x,u,-p)=f(x,u,p).
\end{equation}
In the spirit of the old flow embedding result \cite{SandstedeFiedler} we show that any planar flow
\begin{align}\label{eq(4.2)}
\nonumber\dot{a}&= g(a,b)\\[-2mm]
\hspace*{0,01mm}\\[-2mm]
\nonumber\dot{b}&=h(a,b)
\end{align}
can be realized in this class of PDEs, by embedding (\ref{eq(4.6)}) below, provided that (\ref{eq(4.2)}) is also reflection symmetric, i.e.
\begin{align}\label{eq(4.3)}
\nonumber g(a,-b)&=g(a,b)\\[-2mm]
\hspace*{0,01mm}\\[-2mm]
\nonumber h(a,-b)&=-h(a,b)
\end{align}
Since  there exist reflection symmetric planar vector fields with nonstationary  periodic orbits, PDEs (\ref{eq(1.1)}) with the associated nonlinearity $f$ do not possess Lyapunov functions of the form (\ref{eq(1.4)}), (\ref{eq(1.9)}), (\ref{eq(1.10)}).

\ppar

Our realization of the ODE flow (\ref{eq(4.2)}) will be in the invariant subspace
\[E=\text{span}\{\bc,\bs\}\]
of the first Fourier modes $\bc=\cos x$, $\bs=\sin x$. In fact we define
\begin{equation}\label{eq(4.5)}
f(x,a\bc+b\bs,-a\bs+b\bc):=(a+g(a,b))\bc+(b+h(a,b))\bs
\end{equation}
for all $x\in S^1$ and $a,b\in\R$. This is easy, inverting the rotation
\begin{equation}\label{eq(4.6)}
\begin{pmatrix}u\\ p\end{pmatrix} =\begin{pmatrix}\bc&\bs\\ -\bs&\bc\end{pmatrix}\begin{pmatrix} a\\ b\end{pmatrix}
\end{equation}
and defining, more explicitly,
\begin{align}\label{eq(4.7)}
\nonumber f(x,u,p):=&(u\bc-p\bs+g(u\bc-p\bs,u\bs+p\bc))\bc+\\[-2mm]
\hspace*{0,01mm}\\[-2mm]
\nonumber&+(u\bs+p\bc+h(u\bc-p\bs,u\bs+p\bc))s
\end{align}
for all $x,u,p$. Note how reflection symmetry (\ref{eq(4.3)}) of $g,h$ implies reflection symmetry (\ref{eq(4.1)}) of $f$. Plugging the Ansatz
\begin{equation}\label{eq(4.8)}
u(t,x):= a(t)\bc+b(t)\bs\in E
\end{equation}
into PDE (\ref{eq(1.1)}) we see from (\ref{eq(4.5)}) that such $u$ solve (\ref{eq(1.1)}) iff the coefficients $a(t)$, $b(t)$ satisfy the planar ODE (\ref{eq(4.2)}). This proves our realization claim and establishes the possibility of nonstationary periodic orbits.

\ppar

An analogous construction based on the span of $\cos(nx)$, $\sin(nx)$, instead, shows the possibility of nonstationary periodic orbits in the presence of any finite number of reflection symmetries of PDE (\ref{eq(1.1)}) with respect to $x\in S^1$.

\section{Concluding remarks}\label{sec5}
We briefly comment on the related problem of global attractors for PDE (\ref{eq(1.1)}), on generalizations to quasilinear and nonlinear equations, on finite time blow-up and, finally, on the hidden extension to imaginary $p=u_x$ in our construction of the Lyapunov function $V=\int_0^1 L(x,u,u_x)dx$. 

\ppar

One purpose of Lyapunov functions is to reveal the gradient flow variational character of PDE (1.1) on the circle. In particular we prove convergence to equilibria. Under a dissipativeness assumption on $f$, the \textit{global attractor} $\cA_f$, alias the bounded set of solutions which exist and stay uniformly bounded for all positive and negative times, has received much attention. In presence of a Lyapunov function (\ref{eq(1.4)}), (\ref{eq(1.9)}), (\ref{eq(1.10)}), the global attractor consists of equilibria  and their heteroclinic orbits, only. In contrast, consider the SO(2)-equivariant case $f=f(u,u_x)$ on the circle $x\in S^1$, which does not admit a Lyapunov function. The global attractor $\cA_f$ in this case consists of equilibria, rotating waves, and the heteroclinic orbits connecting them; see \cite{MatanoNakamura}. In \cite{FiedlerRochaWolfrum} the heteroclinic connections were studied by, first, freezing all rotating waves to become circles of nonhomogeneous equilibria and, second, symmetrizing $f$ to become even in $p=u_x$, by suitable homotopies. The present paper then provides an explicit Lyapunov function to deal with the symmetrized case of frozen waves. The main tool in \cite{FiedlerRochaWolfrum} to study the remaining heteroclinic orbits between equilibria was a Sturm nodal property going back to Sturm \cite{Sturm} (1836).  See also \cite{Angenent} and the references there. Hence we call such global attractors $\cA_f$
\textit{Sturm attractors}.

\ppar

Matano in fact studies \textit{quasilinear parabolic PDEs} of the form
\begin{equation}\label{eq(5.1)}
u_t=a(x,u,u_x)u_{xx}+f(x,u,u_x)
\end{equation}
in \cite{Matano}, where $a$ is assumed uniformly positive. The derivation (\ref{eq(2.4)}) -- (\ref{eq(2.11)}) then remains valid if we replace the substitution $u_{xx}=u_t-f$ by $u_{xx}=a^{-1}u_t-a^{-1}f$  there. In particular
\begin{equation}\label{eq(5.2)}
\dot{V}=-\int^1_0a^{-1}L_{pp}(u_t)^2dx
\end{equation}
and we just have to replace $f$ by $f/a$ in (\ref{eq(2.5)}) -- (\ref{eq(2.11)}). Similarly theorem 1.1 remains valid for $O(2)$-equivariant
\begin{equation}\label{eq(5.3)} 
a=a(u,p)=\oa(u,\ffrac p^2)
\end{equation}
if we replace $\of$ by $\of/\oa$ in (\ref{eq(1.18)}) -- (\ref{eq(1.21)}) and replace $L_{pp}$ in (\ref{eq(1.23)}), (\ref{eq(1.24)}) by $L_{pp}/a$.

\ppar

For fully nonlinear parabolic equations
\begin{equation}\label{eq(5.4)} 
u_t=f(x,u,u_x,u_{xx})
\end{equation}
and their equivariant variants a Lyapunov function is not known. Under separated boundary conditions convergence of bounded solutions to single equilibria may still be possible to prove, based on Sturm nodal properties. Albeit the technical ingredients are not sufficiently developed, at present, to provide a short proof here.

\ppar

Returning to the $O(2)$-equivariant semilinear case (\ref{eq(1.1)}) with $f=\of(u,\frac{1}{2}u_x^2)$ it may be interesting to explore the consequences of our Lyapunov function for \textit{blow-up} on the circle $x\in S^1$; see also \cite{FiedlerMatano} for the case of separated boundary conditions. Basically two different effects may occur. First, the Lagrangian integrand $L$ of the Lyapunov function may become unboundedly negative via $u$, in (\ref{eq(1.22)}). Second, the characteristics $q=q(u)$ in (\ref{eq(1.18)}), (\ref{eq(1.19)}) may already explode for finite values of $u$, terminating our very definition of the Lagrangian integrand $L$ in (\ref{eq(1.20)}), (\ref{eq(1.21)}). It will be of interest to compare this second phenomenon, which may occur for nonlinearities $f$ which grow superquadratically in the gradient $u_x$, with the gradient blow-up described in \cite{Kruzhkov}.

\ppar

We conclude with a \textit{complex curiosity} in our construction of the Lagrangian integrand $L$ via the characteristics (\ref{eq(1.18)}), (\ref{eq(1.19)}). Let us first interpret the characteristic
\begin{equation}\label{eq(5.5)}
\frac{d}{du}q=-\of(u,q).
\end{equation}
As long as $q$ remains positive it is easy to see that $q=q(u)$ solves (\ref{eq(5.5)}) iff any solution of
\begin{equation}\label{eq(5.6)}
u_x=\pm\sqrt{2q(u(x))}
\end{equation}
with $u(x)$ in that positivity domain of $q$ solves the equilibrium ODE
\begin{equation}\label{eq(5.7)}
0=u_{xx}+\of(u,\ffrac u^2_x)
\end{equation}
of the PDE (\ref{eq(1.1)}). For a proof just multiply (\ref{eq(5.7)}) by $u_x$ and compare with (\ref{eq(5.5)}) via the chain rule applied to $\frac{d}{dx}q(u(x))$:
\begin{align}\label{eq(5.8)}
\nonumber0&=\frac{d}{dx}(\ffrac u^2_x)+\of(u,\ffrac u^2_x)u_x \quad\text{versus}\\[-2mm]
\hspace*{0,01mm}\\[-2mm]
\nonumber0&=\frac{d}{dx}q(u)\ ~ +\of(u,q(u))u_x\,.
\end{align}
A trivial example, again, are $\of=f(u)$ independent of $q$, where
\begin{equation}\label{eq(5.9)}
q=-F(u)+E
\end{equation}
with the  primitive $F$ of $f$ and the energy $E$ of the second order pendulum equation $u_{xx}+f(u)=0$. Indeed  (\ref{eq(5.6)}) integrates that pendulum, reading
\begin{equation}\label{eq(5.10)}
u_x=\pm\sqrt{2(E-F(u))}.
\end{equation}

\ppar

Our evolution $\Psi^{u_1,u_0}$ of the characteristic equation in  (\ref{eq(1.18)}),  (\ref{eq(1.19)}), however, does not stop at $q=0$. 
Instead it happily proceeds through negative $q=\ffrac p^2$, alias imaginary $p=u_x$, to re-emerge as positive in other regions of the phase plane $(u,q)$. It may therefore become a fascinating speculation to ponder the significance of our simple Lyapunov function for extensions to complex, rather than just real, values of $u$ and $u_x$.

\end{document}